\documentclass[a4paper]{amsart}
\usepackage{amsthm}
\usepackage{xcolor}
\def\P{{\mathbf P}}
\def\Q{{\mathbf Q}}
\def\EE{{\mathfrak E}}
\def\EEE{{\mathcal E}}
\newcommand{\Her}{{\mathcal H}}
\def\FF{{\mathfrak F}}
\def\FFF{{\mathcal F}}

\def\|{\, |\, }
\def\d{{\text{ d}}}
\def\R{{\mathbf R}}
\def\N{{\mathbf N}}
\def\esp#1#2{{\mathbf E}_{#1}\left[{#2}\right]}
\def\Id{{\text{Id}}}
\def\car{{\mathbf 1}}
\def\xdif{\text{ d}}
\newcommand\dist{\operatorname{dist}}
\newcommand{\Lip}{\operatorname{Lip}}
\definecolor{redbox}{RGB}{191,18,56} %
\definecolor{blackbox}{RGB}{0,0,0} %
\definecolor{brownbox}{RGB}{128,99,90} %
\definecolor{macouleur1}{HTML}{D1A42E}
\definecolor{macouleur2}{HTML}{9D7E3C}
\definecolor{macouleur3}{HTML}{9B5C1C}
\definecolor{macouleur4}{HTML}{F9C258}
\usepackage{tikz}
\usetikzlibrary{decorations.pathreplacing,calc}
\usetikzlibrary{arrows,shapes,decorations,automata,backgrounds,petri,shadows,fpu,mindmap}
\newtheorem{thrm}{Theorem}
\newtheorem{dfntn}{Definition}

\begin{document}

\title{The Stein-Dirichlet-Malliavin method}
\author{L. Decreusefond}\address{Institut Mines-Telecom, Telecom ParisTech, CNRS LTCI, Paris, France}

%
%
\begin{abstract}
The Stein's method is a popular method used to derive 
upper-bounds of distances between probability distributions. It can be viewed, in certain of its formulations, as  an
avatar of the semi-group or of the smart-path method used commonly in
Gaussian analysis. We show how this procedure  can be enriched by
Malliavin calculus leading to a functional approach valid in infinite
dimensional spaces.
 \end{abstract}
%
%
\maketitle

\section{Introduction}
\label{sec:introduction}
Distances between probability or probability metrics is a very old topic since it is rich of a
wide range of applications. As mathematical objects, it is natural to
define a metric topology on spaces of probability measures. As
modeling objects, it is natural to compare probability measures which
appear in the mathematical representations of random phenomena. This
topic has at least three facets: The diverse definitions of
probability metrics which are tailored for each applications; the
computations and comparisons of these different distances for the
widest possible range of situations and at last, the applications
which go from mathematical considerations like functional inequalities
to more practical results of rate of convergence of stochastic algorithms.
The Figure~\ref{fig_decreusefond_esaim:mindmap} shows a partial view
of the different aspects of this subject.

\begin{figure}[!h]
  \centering
\begin{tikzpicture}
\path [small mindmap, text=white,concept color=macouleur1]
    node[concept ] {Probability\\ metrics}
    [clockwise from=0]
   child[concept color=macouleur2]
          {node[concept] {Computations}
 [clockwise from=22]
            child[concept color=blue!60]
               {node[concept] {Girsanov\\ $\Q\ll \P$}}
            child[concept color=red]
               {node[concept] {Stein\\ $\Q=T^*m$}}
}
   child[concept color=macouleur2] 
         {node[concept] {Types}
              child[concept color=macouleur3]
                   {node[concept] {Optimal\\ Transport}
                      [clockwise from=-10]
                   child[concept color=red]
                      {node[concept] {Rubinstein}}
                   child[concept color=blue!60]
                      {node[concept] {Wasserstein}}
}
              child[concept color=macouleur3]     
                   {node[concept] {Prohorov}}
              child[concept color=blue!60]     
                   {node[concept] {Entropy}}     
                 }
    child[concept color=macouleur2]
       {node[concept] {Applications}
  [clockwise from=-45]
         child[concept color=blue!60]
         {node[concept] {Functional \\ inequality}}
         child[concept color=red]
         {node[concept] {Convergence \\ rate}}
         child[concept color=macouleur3]
         {node[concept] {Ergodicity}}
}  ;
\end{tikzpicture}
  \caption{Mindmap}
  \label{fig_decreusefond_esaim:mindmap}
\end{figure}

A few words are in order to explain the blue and red colors. For the
computations of distances between measures $\mu$ and $\nu$, we need to impose some
relationships between these two measures. Absolute continuity is one
very frequent type of relationships between two measures. The
Radon-Nykodim theorem gives a precious tool to estimate
divergence-like and Wasserstein distances (see for instance
\cite{MR2036490} for such an application). One may also reverse the
point of view: Given a positive function $F$, compare the $\mu$ and
$\nu=F\d\mu$ to obtain some precious functional inequalities on $F$
(see \cite{MR3155209}). These results thus belong to the same
\textsl{spirit} and are colored in blue. Another natural way to put a
structure between two measures is to have a map which transforms a
known measure into another one and to compare this transformed measure to
a reference probability. This is exactly the framework in which the
Stein's method performs well if we consider Kantorovitch-Rubinstein
type distances (defined below). Typical applications of these form of
distances are to give the convergence rates of celebrated theorem like
CLT or Berry-Esseen Theorem or of random algorithms
\cite{MR99k:28007}. The links between these different points justify
that they are all  colored in red. 

This paper is a rather informal introduction to the
Stein-Dirichlet-Malliavin method (SDM for short henceforth). This is an extension of the
classical Stein's method, enriched by the structure given by Dirichlet
forms and Malliavin calculus. We hope that this new point of view will
lead to more systematic proofs of convergence, extending  their
applicability. The price to pay is to master some new concepts from
Malliavin calculus like the gradient and its associated
adjoint. That is why we tried to maintain the technicalities at the
lowest possible level, insisting more on the ideas at play.

We first show the different kinds of probability metrics that exist in
the literature. We do not pretend to be exhaustive but aim to point
out to the wide diversity of possible definitions. In Section
\ref{sec:probability-metrics}, we establish the principles of the SDM
method and show how it can be applied to the Poisson-Gaussian
convergence. We then explain how to construct the necessary structures
to extend this procedure to infinite dimensional spaces. In
Section~\ref{sec:edgeworth-expansion}, Edgeworth expansions are
obtained by iterating the previous procedure as often as desired.

\section{Taxonomy of probability metrics}
\label{sec:probability-metrics}
In what follows, all the probability measures are defined on
Polish spaces denoted either by $\EE$ or $\FF$, whose
borelian $\sigma$-fields is $\mathfrak B(\EE)$, respectively $\mathfrak B(\FF)$. 
There are several notions of metrics between probability
measures. An interesting survey of the main variants
and their mutual relationships can be found in \cite{INSR:INSR419}. Each of one is often adapted to a particular purpose. They
can roughly and partly be classified in three types. The first one is
the so-called Prokhorov distance.
\begin{equation*}
 \text{Dist}_{\text{Pro}}(\P,\Q)=\inf\Bigl\{\epsilon >0, 
\P(A)\le \Q(A^\epsilon)+\epsilon \text{ for all } A\in \mathfrak B(\EE)\Bigr\},
\end{equation*}
where $A^\epsilon$ is the $\epsilon$-neighborhood of $A$ defined by
\begin{math}
  A^\epsilon=\{y\in \EE,\ \exists x\in A, d(x,y)\le \epsilon\}.
\end{math}
This distance is crucial as its associated topology is precisely the topology of the
convergence in distribution, i.e. we have the following theorem which
can be found in \cite{MR982264}.
\begin{thrm}
  A sequence $(\P_n,\, n\ge 1)$ of probability measures converges
  weakly to $\P$ if and only if $ \text{Dist}_{\text{Pro}}(\P_n,\P)$
  tends to $0$ as $n$ goes to $\infty$.
\end{thrm}
Unfortunately, this distance is hardly computable and that justifies
the search for alternative and more tractable definitions.
A vast category of probability metrics is represented by the
$f$-divergence defined as follows.
\begin{dfntn}
    Let  $f$ be a convex function such that $f(1)=0$. Then, for two
    probability measures $\P$ and $\Q$ on a Polish space $\EE$, 
\begin{equation*}
D_f(\Q \| \P)=
\begin{cases}
\displaystyle{\int_{\EE}} f\left(\frac{\xdif\Q}{\xdif\P}\right)\xdif \P      & \text{ if } \Q\ll
\P,\\
\infty & \text{ otherwise.}  
\end{cases}
    \end{equation*}
\end{dfntn}
For instance, if we choose $f=t\ln t$, we obtain the Kullblack-Leibler
distance. The Hellinger distance corresponds to the case where
$f(t)=(\sqrt{t}-1)^2$. Total variation between absolutely continuous
measures boils down to take $f(t)=|t-1|$.

Another class of distances between measures can be obtained via
optimal transportation theory. For general results about this theory,
we refer to the books \cite{MR99k:28006,MR99k:28007,Villani:2007fk,MR1964483}.
\begin{dfntn}
  Let $(\EE, \, \P)$ and $(\FF,\, \Q)$ two Polish spaces equipped with
  a probability measure and $c$ a semi-continuous
  function from $\EE\times \FF$ to $\R^+\cup \{\infty\}$. The
  optimal-transportation problem or Monge-Kantorovitch problem
  $\text{MKP}(\P,\Q)$ is to find
  \begin{equation*}
    \min_{\gamma \in \Sigma (\P,\Q)}\int_{\EE\times \FF} c(x,y)\xdif \gamma(x,y)
  \end{equation*}
where $ \Sigma (\P,\Q)$ denoted the space of probability measures on
$\EE\times \FF$ with first marginal $\P$  and second marginal $\Q$.
\end{dfntn}
Said otherwise in a more probabilistic way, it amounts to find the
coupling between $\P$ and $\Q$ which minimizes the cost, i.e. to
construct on the same probability space, two random variables $X$ and
$Y$ of respective distribution $\P$ and $\Q$ which minimizes
$\esp{}{c(X,Y)}$ among all the possible constructions. The usual cost
functions are of the type $c(x,y)=\text{dist}(x,y)^p$ where
$\text{dist}$ is a distance and $p$ a positive real number. For the
Euclidean distance and $p=2$, we can construct the so-called
Wasserstein  distance by considering
\begin{equation*}
  W(\P,\Q)=\sqrt{\min_{\gamma \in \Sigma (\P,\Q)}\int_{\R^d\times \R^d} \|x-y\|^2\xdif \gamma(x,y)}.
\end{equation*}
All the distances viewed so far are not unrelated as many functional
inequalities do exist between all of them. Just to mention two examples, the
Pinsker inequality  states that the total variation distance is
controlled by the Kullblack-Leibler distance.
\begin{equation*}
  D_{|t-1|}(\P,\Q)\le \sqrt{\frac12 D_{t\ln t}(\P,\Q) }.
\end{equation*}
On the other hand, the so-called HWI identity (see \cite{MR1964483}) relates the relative
entropy (H), the Wasserstein distance (W) and the Fischer information
(I) as follows. 
\begin{thrm}
  \label{sec:scope-steins-method-1} Let $\P$ and $\Q$ two probability
  measures on $\R^n$ such that $\P=\exp(-V)\xdif x$ with $\nabla^2V\ge
  K\Id_n$. Then,
\begin{equation*}
      D_{t\ln t}(\P,\Q)\le W(\P,\Q)\sqrt{D_{\nabla |\ln t|^2}(\P,\Q)}-\frac{K}{2} W(\P,\Q)^2.
    \end{equation*}
\end{thrm}
These examples are here only to give a glimpse of  the vast subject
of the relationship between all these notions of distances. However,
this is not the true subject of the present paper. 
The theorem which justifies the sequel is known as
Kantorovitch-Rubinstein theorem (see \cite{MR982264,MR622552}) and
says the following.
\begin{thrm}
  For $\P$ and $\Q$ two probability measures on a Polish space $\EE$,
   consider the Monge-Kantorovitch problem for a cost function $c$ 
  which is a distance on $\EE$. Then, we have  the following representation
  \begin{equation*}
\min_{\gamma \in \Sigma (\P,\Q)}\int_{\EE\times \FF} c(x,\, y)\xdif \gamma(x,y)=\sup_{F \in \Lip_c(1)}\left(\esp{\P}{F}-\esp{\Q}{F}\right),
  \end{equation*}
where $F \in \Lip_c(1)$ means that $F$ is $c$-Lipschitz continuous:
$|F(x)-F(y)|\le c(x,y)$ for all $x,\, y\in E$. The resulting distance
between $\P$ and $\Q$, 
will be called henceforth the Kantorovitch-Rubinstein distance as in \cite{MR1964483}.
\end{thrm}
This formulation of a distance motivates alternative definitions by
changing the set of test functions. For instance, for
$\FFF=\{\car_{(-\infty;\, x]},\, x\in R\}$,
  \begin{equation*}
    \sup_{F \in \FFF}\left|\esp{\P}{F}-\esp{\Q}{F}\right|
  \end{equation*}
is the total-variation distance. It turns out that Stein's method is
particularly well suited to estimate such kind of distances as we
shall see now.


\section{Stein's method}
\label{sec:steins-method}
Historically, the Stein's method for Gaussian distribution dates back
to the seminal paper of Stein \cite{stein1972}. It was soon extended
to the Poisson distribution in the paper of Chen \cite{MR0370693}. It
is then impossible to track all the extensions of this approach, made
mainly by A. Barbour and his collaborators, to several other
distributions like compound Poisson \cite{MR1920275}, Poisson point processes\cite{MR2001a:60058}, stationary measure of birth-death
process, even Brownian motion \cite{MR1035659}. For a whole account of
all this period, one may refer to the books
\cite{MR1708412,MR93g:60043} and references therein. The main
breakthrough came with the paper of Nourdin and Peccati
\cite{MR2520122}, in which it is shown that combining Malliavin
calculus and Stein's approach, one can obtain a rather simple proof of
the striking \textsl{fourth moment theorem}, established earlier in
\cite{MR2118863}. This was the starting point of a bunch of articles
with with a wide area of applications: rate of convergence in the
central limit theorem, Berry-Esseen theorem, iterated-logarithm
theorem, limit theorems on manifolds, etc.


\subsection{Dirichlet-Malliavin structure}
\label{sec:going-functional}
The procedure of the Stein's method can be abstracted within the setting of Dirichlet
structures (for details, we refer to
\cite{bouleau-hirsch,MR569058,MR1214375}). The subsequent explanations
are at a very formal level since the hard part for this machinery to work
is to find the convenient functional spaces for each case of
applications.

The first idea underlying the Stein's method is to characterize the target
measure by an algebraic equation: Find a functional operator $L$ on
$\mathcal F$ such that $\esp{\Q}{LF}=0$ for any $F$ in $\mathcal F$ if
and only if $\Q=\P$.
It turns out that this functional operator $L$  can be
viewed as the (infinitesimal) generator of a Markovian semi-group,
which we denote by $P=(P_t,\, t\ge 0)$ whose stationary measure is $\P$:
The image measure of $\P$ by $P_t$ is still $\P$ for any $t\ge
0$. Under some technical hypothesis, there exists a strong ergodic Markov
process $X=(X(t),\, t\ge 0)$ of invariant measure $\P$ and of generator
$L$. It must be noted that the knowledge of one of $L$, $P$ or $X$ is
equivalent to the knowledge of the other two. Formally speaking, for
any $x\in E$,
\begin{equation*}
  P_tf(x)=e^{tL}f(x), \ Lf(x)=\left.\frac{dP_tf(x)}{dt}\right|_{t=0}, \ P_tf(x)=\esp{}{f(X(t))\, |\, X(0)=x}. 
\end{equation*}
One can also associate to $X$, the so-called Dirichlet form defined
formally by 
\begin{equation*}
  \EEE(F,G)=\esp{\P}{LF\, G},
\end{equation*}
for any $F$ and $G$ sufficiently regular. As before, if we are given
such a bilinear form $\EEE$, one can retrieve $L$ by the following relationship: For
any $F$, $LF$ is the unique element $H$  such that for any
$G$, $\EEE(F,\, G)=\esp{\P}{HG}$. This means that whichever of $L$,
$X$, $P$ or $\EEE$ we are given, the others are uniquely
determined (the reader is referred to the particularly illuminating 
Diagram 2, page 36 of \cite{MR1214375}). 

 Within this framework, it
is easy to see that the Stein-Dirichlet representation formula holds:
For any bounded $F$,
\begin{equation}
\label{eq_decreusefond_esaim:1}
  \esp{\Q}{F}-\esp{\P}{F}=\esp{\Q}{\int_0^\infty LP_tF \d t}.
\end{equation}
This formula is also known as \textsl{the semi-group method} or
\textsl{the smart-path formula} in the Stein's method literature.
This means that we can write
\begin{equation*}
  \dist_{\mathcal F}(\P,\Q)=\sup_{F\in \mathcal F}\left| \esp{\P}{F}-\esp{\Q}{F}\right|=\sup_{F\in \mathcal F}\left|\esp{\Q}{\int_0^\infty LP_tF \d t}
\right|.
\end{equation*}
Instead of using coupling arguments to estimate this right-hand-side
as usually done in the Stein's method,
we use another functional operator which is the 
\textsl{gradient} in the sense of Malliavin. It is  usually denoted by
$D$ and satisfies the identity  $L=D^*D$ where
$D^*$ is the adjoint of $D$. This a square root of the symmetric
operator $L$, but not all square-roots are interesting as we also  need 
a nice commutation relationship between $D$ and $P$. A few examples
are the best way to illustrate what we mean.

\subsection{One dimensional examples}
\label{sec:gaussian-measures}

If $\P$ denote the standard Gaussian measure on $\R$, then $X$ is the
Ornstein-Uhlenbeck process defined by 
\begin{equation*}
  dX(t)=\sqrt{2}\d B(t)-X(t)\d t, \ X(0)=x,
\end{equation*}
where $B$ is a standard one-dimensional Brownian motion. 
A  straightforward application of the It\^o formula gives the following
expression of $X$:
\begin{equation*}
  X(t)=e^{-t}x+\sqrt{2}\int_0^te^{-(t-s)}\d B(s).
\end{equation*}
It is then easy to see that $X(t)\sim\mathcal N(e^{-t}x,\, 1-e^{-2t})$, which, in turn,
entails the Mehler
representation formula:
\begin{equation*}
  P_tF(x)=\int_\R F(e^{-t}x+\sqrt{1-e^{-2t}}y)\d\P(y).
\end{equation*}
It follows by differentiation and integration by parts that for $F\in
\mathcal C^2_b$,
\begin{equation*}
  LF(x)=xF^\prime(x)-F^{\prime\prime}(x), \text{ for all } x\in \R.
\end{equation*}
The \textsl{Malliavin} gradient is the usual derivative operator and
standard computations show that 
\begin{equation*}
  \int_\R DF(x) \, G(x)\d\P(x)=\int_\R F(x)(xG(x)-DG(x))\d\P(x),
\end{equation*}
hence that $D^*G(x)=xG(x)-DG(x)$ and $L=D^*D$. Moreover, we have
$DP_tF(x)=e^{-t}P_tDF(x)$ which is the commutation relationship
alluded above.

If $\P$ represents the Poisson measure on $\N$ of parameter $\lambda$, the process $X$ can be
viewed as the number of occupied servers in an M/M/$\infty$ queue (see
\cite{Decreusefond:2012sys}), $L$ is the corresponding generator:
\begin{equation*}
  LF(x)=\lambda(F(x+1)-F(x))+x(F(x-1)-F(x)), \text{ for all } x\in \N,
\end{equation*}
with the convention that $0.F(-1)=0$.
The \textsl{gradient} is defined by 
\begin{equation*}
  DF(x)=F(x+1)-F(x),
\end{equation*}
and  we have $DP_tF=e^{-t}P_tDF$ (see \cite[Theorem
11.16]{Decreusefond:2012sys} or \cite{DST:functional}). 
For the scalar product in $L^2(\P)$, we have
\begin{equation}\label{eq_decreusefond_esaim:3}
  \int_\N DF(x)\, G(x)\d\P(x)=\int_\N
  F(x)(\frac{x}{\lambda}G(x-1)-G(x))\d \P(x).
\end{equation}
Hence, 
\begin{equation*}
 D^*F(x)=\frac{x}{\lambda}G(x-1)-G(x) \text{ and } L=D^*D.
\end{equation*}
We now show how these constructions do articulate  to give a new
approach to the Stein's method.

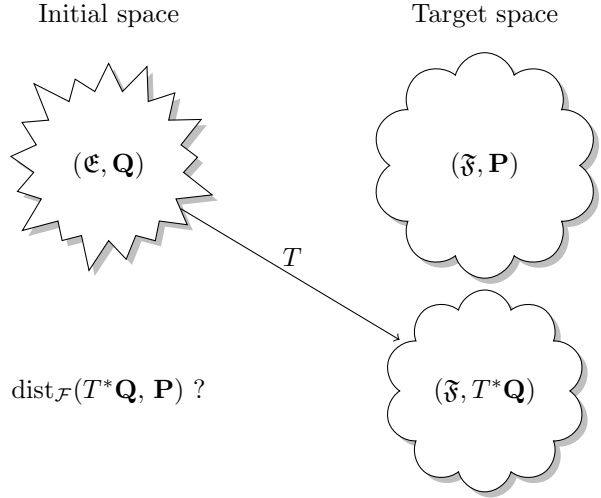
\begin{figure}[!h]
  \centering
   \begin{tikzpicture}
\node at (0,2) {Initial space};
\node at (5,2) {Target space};
    \node[starburst,minimum size=3cm,drop shadow,fill=white,draw] (a) at (0,0) {$(\EE,\Q)$};
    \node[cloud,minimum size=3cm,drop shadow,fill=white,draw] (b) at (5,0) {$(\FF,\P)$};
{\node[cloud,drop shadow,fill=white,draw] (c) at (5,-3) {$(\mathfrak F,T^*\Q)$};
\path[->,shorten >=2pt] (a) edge node[above] {$T$} (c);}
{\node at (0,-3) {$\text{dist}_{\mathcal F} (T^*\Q,\, \P)$ ?};}
  \end{tikzpicture}    
  \caption{Comparison between a measure $\P$ and $T^*\Q$.}
  \label{fig_decreusefond_esaim:comparaison}
\end{figure}

It is well known that for $Z_\lambda$ a Poisson random variable of
parameter $\lambda$, 
\begin{equation*}
  \hat{Z}_\lambda
  =\frac{Z_\lambda-\lambda}{\sqrt{\lambda}}\xrightarrow{\lambda \to
    \infty} \mathcal N(0,\ 1) \text{ in distribution.}
\end{equation*}
We are going to use the Stein-Dirichlet-Malliavin method to evaluate
the rate of convergence. We are in a situation where the target measure in defined $\R$ whereas
the initial randomness comes from a probability measure on $\N$. The
map $T$ defined by 
\begin{align*}
  T\, :\, \EE=\N&\longrightarrow \FF=\R\\
n&\longmapsto \frac{n-\lambda}{\sqrt{\lambda}},
\end{align*}
maps one space to the other and we are to evaluate the distance
between $T^*\Q_\lambda$, the image measure of $\Q_\lambda$, the Poisson($\lambda$)
probability, by the map $T$ and $\P$ the standard normal distribution
on $\R$. This is a particular case of the  general situation
illustrated in Figure~\ref{fig_decreusefond_esaim:comparaison}.

In view of \eqref{eq_decreusefond_esaim:1}, we have to estimate
\begin{equation*}
  \sup_{F\in \mathcal F} \int_0^\infty
\int_\R  x.(P_tF)^\prime(x)-(P_tF)^{\prime\prime}(x)
\d T^*\Q_\lambda(x) \d t,
\end{equation*}
where is the Ornstein-Uhlenbeck semi-group given by the Mehler formula
above and $\mathcal F$ is a functional space to be conveniently chosen. According to the definition of $T$, the quantity to maximize is
equal to 
\begin{equation*}
  \esp{}{\int_0^\infty
    \hat{Z}_\lambda.(P_tF)^\prime(\hat{Z}_\lambda)-(P_tF)^{\prime\prime}(\hat{Z}_\lambda)\d
  t}.
\end{equation*}
Applying \eqref{eq_decreusefond_esaim:3} to $G=1$ and $ F\circ T$, we
get
\begin{equation*}
\sqrt{\lambda}  \ \esp{}{F(\hat{Z}_\lambda+\frac{1}{\sqrt{\lambda}})-F(\hat{Z}_\lambda)}=\esp{}{\hat{Z}_\lambda\,F(\hat{Z}_\lambda)}.
\end{equation*}
Hence, 
\begin{equation}\label{eq_decreusefond_esaim:5}
  \esp{}{\hat{Z}_\lambda.(P_tF)^\prime(\hat{Z}_\lambda)}=\sqrt{\lambda}  \esp{}{(P_tF)^\prime(\hat{Z}_\lambda+\frac{1}{\sqrt{\lambda}})-(P_tF)^\prime(\hat{Z}_\lambda)}.
\end{equation}
For any $t>0$, the regularizing properties of $P_t$ entails that
$P_tF$ is thrice differentiable. Hence, 
\begin{equation}\label{eq_decreusefond_esaim:6}
 (P_tF)^\prime(\hat{Z}_\lambda+\frac{1}{\sqrt{\lambda}})-(P_tF)^\prime(\hat{Z}_\lambda)=
 \frac{1}{\sqrt{\lambda}} (P_tF)^{\prime\prime}(\hat{Z}_\lambda)+\frac{1}{{\lambda}} \int_0^1 (1-r)
 (P_tF)^{(3)}(\hat{Z}_\lambda+\frac{r}{\sqrt{\lambda}})\d r.
\end{equation}
And then, a miracle occurs: The term involving the second order
derivative vanishes and we are lead to maximize 
\begin{equation}\label{eq_decreusefond_esaim:4}
 \frac{1}{\sqrt{\lambda}}\  \esp{}{\int_0^\infty  \int_0^1 (1-r)
 (P_tF)^{(3)}(\hat{Z}_\lambda+\frac{r}{\sqrt{\lambda}})\d r \d t}
\end{equation}
for $F$ over $\mathcal F$. There is now a delicate point. If $F$ is in
$\mathcal C^1_b$, we already mentioned that 
\begin{equation*}
  (P_tF)^\prime(x)=e^{-t}P_t(F^\prime)(x).
\end{equation*}
Furthermore, by integration by parts with respect to the Gaussian
measure, it is easy to see that
\begin{equation*}
  (P_tF)^{(k)}(x)=\left(\frac{e^{-t}}{\sqrt{1-e^{-2t}}}\right)^k\int_\R
  F(e^{-t}x+\sqrt{1-e^{-2t}}y)\, y^k\d\P(y),
\end{equation*}
whenever $F$ is bounded, for any $k\ge 1$. At first glance, it seems
easy to bound \eqref{eq_decreusefond_esaim:4} by using the previous
formula for $k=3$. Unfortunately, the term $\exp(-kt)(1-\exp(-2t))^{-k/2}$ is
integrable over $[0,+\infty)$ only for $k=1$. Hence, we must choose
$\mathcal F=\{F\in \mathcal C^2_b, \Vert F\Vert_{\mathcal C^2_b}\le 1\}$ and then we have
\begin{multline*}
   \left|(P_tF)^{(3)}(x)\right|=\left|\frac{e^{-3t}}{\sqrt{1-e^{-2t}}} \int_\R
  F^{(2)}(e^{-t}x+\sqrt{1-e^{-2t}}y)\, y\d\P(y)\right|\\ \le
\frac{e^{-3t}}{\sqrt{1-e^{-2t}}}\Vert F^{(2)}\Vert _\infty\ \int_\R |y|\d\P(y).
\end{multline*}
Plugging this inequality into \eqref{eq_decreusefond_esaim:4}, we get 
\begin{multline}\label{eq_decreusefond_esaim:8}
  \sup_{\Vert F\Vert_{\mathcal C^2_b}\le 1}\left|
  \esp{}{F(\hat{Z}_\lambda)}-\int F\d\P\right|
\\ \le  \frac{1}{\sqrt{\lambda}}\int_0^1 (1-r)\d r\ 
  \int_0^\infty\frac{e^{-3t}}{\sqrt{1-e^{-2t}}}\ dt \ \int_\R
  |y|\d\P(y)
=\frac{\sqrt{\pi}}{4\sqrt{2}}\ \frac{1}{\sqrt{\lambda}}\cdotp
\end{multline}
Hence we have established the rate of convergence for the
Kantorovitch-Rubinstein distance associated to $\mathcal F=\{F\in
\mathcal C^2_b, \Vert F\Vert_{\mathcal C^2_b}\le 1\}$. In dimension
$1$, for Gaussian approximation, we could have used
$LF(x)=xF(x)-F^\prime(x)$ as a characterizing operator and thus used
only $1$-Lipschitz functions with a slightly different constant in
front of the $\lambda^{-1}$ factor, namely
\begin{equation*}
  \sup_{F\in \Lip(1)}\left|
  \esp{}{F(\hat{Z}_\lambda)}-\int F\d\P\right|\le \frac{1}{\sqrt{2\pi}}\ \frac{1}{\lambda}\cdotp
\end{equation*}
Note that this upper-bound is better than the bound obtained by the
classical Stein's method where $(2\pi)^{-1/2}$ is replaced by $1$.
 However, this line of thought is not
applicable to higher dimensions.

More generally, the recipe of the Stein-Dirichlet-Malliavin method is
the following.
\begin{itemize}
\item Characterize the target measure as the stationary distribution
  of an ergodic Markov process,
\item Construct the two Dirichlet-Malliavin structure on both initial
  and target spaces,
\item Perform an integration by parts on the initial space (see
  \eqref{eq_decreusefond_esaim:5}),
\item Replace the gradient on the initial space by a function of the
  gradient on the target space (this is done here by the Taylor
  formula \eqref{eq_decreusefond_esaim:6}), at the price of additional
  terms to be controlled,
\item Finish the computations in the target space using the commuting
  relationship : $DP_t=e^{-t}P_t D$.
\end{itemize}
\subsection{Higher dimensions}
\label{sec:higher-dimensions}
This procedure can be generalized to any dimension provided that we
have Dirichlet-Malliavin structures on both the initial and the target spaces.
For the Gaussian measure in dimension $d$, the generator is given by 
\begin{equation}\label{eq_decreusefond_esaim:2}
   LF(x)=x.D F(x)-\Delta F(x), \text{ for all } x\in \R^d,
\end{equation}
where $D$ is the usual gradient in $\R^d$ and $\Delta$ is the
Laplacian operator. The Mehler formula stays formally the same with an
integral over $\R^d$ instead of $\R$ and  $X$ is the
$\R^d$-valued process composed of $d$ independent copies of the one dimensional
Ornstein-Uhlenbeck process. The Malliavin gradient is still the usual
gradient and the commutation relationship between $D$ and $P_t$ is
easily seen to hold again. We can then retrieve the results of \cite{MR2727319}. 

Real difficulties arise when we try to generalize this approach to
infinite dimensional spaces like the Wiener space. It is tempting to
define $L$ formally as in \eqref{eq_decreusefond_esaim:2}, replacing
the Laplacian by the trace of $D\circ D$. Unfortunately, for this
trace term to exist, we need to restrict the space $\mathcal F$ of
test functions and to choose conveniently the space $\FF$. There are
actually two papers which address this problem. In both of them
\cite{CD:2012,Shih20111236}, despite apparent dissimilarities, we end
by considering $\FF$  a Hilbert space with a Gaussian measure.

Let us show how it works on an example. For $N_\lambda$ a Poisson process on $\R^+$ of intensity
$\lambda$, it is known that 
\begin{equation*}
 \hat{N}_\lambda(t)=\frac{N_\lambda(t)-\lambda
   t}{\sqrt{\lambda}}\xrightarrow{\lambda\to \infty} B(t) \text{ in
   distribution}, 
\end{equation*}
where $B$ is a standard Brownian motion and the convergence is
understood to hold in $\mathbb D$, the Skorohod space of rcll
functions. To compare the two distributions implies to find a common
Hilbert space which supports both the distribution of $B$ and
$\hat{N}_\lambda$. In principle, any Sobolev-like space should do. In
\cite{CD:2012},  we chose the so-called Besov-Liouville space
$I^{\beta, 2}$ for $\beta<1/2$ defined by
\begin{equation*}
  I^{\beta, 2}=\{f, \, \exists \dot f\in L^2([0,1]) \text{ such that }
  f(x)=\frac{1}{\Gamma(\beta)}\, \int_0^x (x-y)^{\beta-1}\dot f(y)\d y\}.
\end{equation*}
It is a Hilbert space when equipped with the scalar-product $\langle
f,\, g\rangle_{\beta,2}=\langle \dot f,\, \dot g\rangle_{L^2}$.
The Wiener measure on this space, denoted by  $\mu_\beta$, is defined by 
\begin{equation*}
  \esp{\mu_\beta}{\exp(i\langle \eta,\, \omega\rangle_{{\beta,\,
        2}})}=\exp(-\frac 12 \langle V_\beta\eta,\,
  \eta\rangle_{{\beta,\, 2}}).
\end{equation*}
where 
\begin{gather*}
I_{0^+}^\beta f(x)=\frac{1}{\Gamma(\beta)}\, \int_0^x
(x-y)^{\beta-1}\dot f(y)\d y,\ I_{1^-}^\beta f(x)=\frac{1}{\Gamma(\beta)}\, \int_x^1
(y-x)^{\beta-1}\dot f(y)\d y\\
\text{ and }   V_\beta=I_{0^+}^\beta\circ I_{0^+}^{1-\beta}\circ
  I_{1^-}^{1-\beta}\circ I_{0^+}^{-\beta}.
\end{gather*}
  The Ornstein-Uhlenbeck semi-group on $( I^{\beta, 2},\, \mu_\beta)$ is
  defined for any $F\in L^2( I^{\beta, 2},\, \mu_\beta)$ by
  \begin{align*}
    P_t^\beta F(u)&:=\int_{I^{\beta, 2}} F(e^{-t}u+\sqrt{1-e^{-2t}}\, v)\d
    \mu_\beta(v).
  \end{align*}
The gradient is the Fr\'echet gradient on $I^{\beta, 2}$ and all the
other properties still holds formally as in finite dimension.

As initial space, we consider $\EE=\mathfrak N$, the space of locally
finite configurations on $\R^+$ equipped with the vague topology. The
measure $\mathbf Q_\lambda$ is such that the canonical process,
denoted by $N_\lambda$, is a
Poisson process of intensity $\lambda$, for details we refer to \cite{CD:2012}.
On the initial space, we actually only need to know the gradient and
an integration by parts formula. Here, we take
\begin{equation*}
  D_xF(N_\lambda)=F(N_\lambda+\delta_x)-F(N_\lambda),
\end{equation*}
where $N_\lambda+\delta_x$ is the configuration $N_\lambda$ with an
additional atom at location $x$. The well-known Campbell-Mecke formula
(\cite{kallenberg83,MR1113698}) is
equivalent to say that 
\begin{equation*}
\esp{\Q_\lambda}{F\ \int_0^1 G_\tau (\d N_\lambda(\tau)-\lambda \d \tau))}=\lambda\ \esp{\Q_\lambda}{\int_0^1
    D_\tau F\  G_\tau\d \tau},
\end{equation*}
for $G$ a deterministic process. The map $T$ is defined by 
\begin{align*}
  T\, :\, \, \mathfrak N& \longrightarrow I^{\beta,2}\\
N & \longmapsto (t\mapsto \frac{N(t)-\lambda t}{\sqrt{\lambda}})\cdotp
\end{align*}
Proceeding exactly along the same
lines as before, one can show that there exists $c_\beta>0$ such that
\begin{equation}\label{eq_decreusefond_esaim:7}
  \sup_{\Vert F\Vert _{ \mathcal C^2_b(I^{\beta,2};\, \R)}\le 1}\left|
    \esp{\Q_\lambda}{F}-\esp{\mu_\beta}{F}\right| \le \frac{c_\beta}{\sqrt{\lambda}},
\end{equation}
where $ \mathcal C^2_b(I^{\beta,2};\, \R)$ is the set of twice
Fr\'echet differentiable functionals on $I^{\beta,2}$, with bounded
differentials. This is the generalization we could expect of \eqref{eq_decreusefond_esaim:8}.

Other examples of the application of this procedure, involving other
functional spaces,  can be found in
the papers \cite{CD:2012,DST:functional}. A similar approach with
Malliavin calculus replaced by a coupling argument appears in \cite{schumacher14}.


\section{Edgeworth expansion}
\label{sec:edgeworth-expansion}
The Stein's method as developed here can be iterated to obtain
Edgeworth expansions. We now want to precise the expansion obtained in
\eqref{eq_decreusefond_esaim:8}. For, we go one step further in
the Taylor formula \eqref{eq_decreusefond_esaim:6}:
\begin{equation*}
  \psi(\hat{Z}_\lambda+1/\sqrt{\lambda})-\psi(\hat{Z}_\lambda)=\frac{1}{\sqrt{\lambda}}\psi^\prime(\hat{Z}_\lambda)+\frac{1}{2\lambda}\psi^{\prime\prime}(\hat{Z}_\lambda)+\frac{1}{6\lambda^{3/2}}\psi^{(3)}(\hat{Z}+\theta/\sqrt{\lambda}).
\end{equation*}
Hence,
\begin{multline}\label{eq_gaussapp20:4}
  \esp{}{ \hat{Z}_\lambda
    DP_tF(\hat{Z}_\lambda)-D^{(2)}P_tF(\hat{Z}_\lambda)}\\ = \frac{1}{2\sqrt{\lambda}}\esp{}{D^{(3)}P_tF(\hat{Z}_\lambda)}+\frac{1}{6\lambda}\esp{}{D^{(4)}P_tF(\hat{Z}+\theta/\sqrt{\lambda})}.
\end{multline}
If $F$ is thrice differentiable with bounded derivatives then $P_tF$
is four times differentiable, hence the last
term of \eqref{eq_gaussapp20:4} is bounded by
$\lambda^{-1}\frac{e^{-4t}}{\sqrt{1-e^{-2t}}}\Vert F^{(3)}\Vert_\infty/6$. Moreover,
applying \eqref{eq_decreusefond_esaim:8} to $DP_tF$ shows that
\begin{equation*}
  \esp{}{D^{(3)}P_tF(\hat{Z}_\lambda)}=
  \esp{\P}{D^{(3)}P_tF}+O(\lambda^{-1/2}).
\end{equation*}
Combining the last two results, we obtain that for $F$ thrice
differentiable
\begin{equation*}
  \esp{}{F(\hat{Z}_\lambda)}-\esp{\P}{F}=\frac{1}{2\sqrt{\lambda}}\esp{\P}{\int_0^\infty
    D^{(2)}P_tF\d t} +O(\lambda^{-1}).
\end{equation*}
This line of thought can be pursued at any order provided that $F$ is
assumed to have sufficient regularity and we get an Edgeworth
expansion up to any power of $\lambda^{-1/2}$.  Using the properties
of Hermite polynomials, this leads to the expansion:
\begin{equation*}
  \esp{}{F(\hat{Z}_\lambda)}-\esp{\P}{F}
  =\frac{1}{6\sqrt{\lambda}}\esp{\P}{F \Her_3}+O(\lambda^{-1}),
\end{equation*}
where $\Her_n$ is the $n$-th Hermite polynomials. In \cite{CD:2014},
we generalized this approach to the Poisson process-Brownian motion
convergence established in \eqref{eq_decreusefond_esaim:7}.
\section{Conclusion}
\label{sec:conclusion}

We showed how the Stein's method can be abstracted in the framework of
Dirichlet forms and Malliavin calculus. This gives raise to a new
method of proof which can be applied to infinite dimensional spaces
and iterated to get Edgeworth expansions. One open question is to apply this approach to other
limiting processes like stable or max-stable processes, Brownian
bridges, etc.

\end{document}